\documentclass{amsart}

\title[Stability of Spacelike Hypersurfaces]{Stability of Spacelike Hypersurfaces in Foliated Spacetimes}

\author{A. Barros, A. Brasil and A. Caminha}

\usepackage{graphicx,pstricks,pst-plot}

\newtheorem{theorem}{Theorem}[section]
\newtheorem{lemma}[theorem]{Lemma}
\newtheorem{proposition}[theorem]{Proposition}
\newtheorem{corollary}[theorem]{Corollary}

\theoremstyle{definition}
\newtheorem{definition}[theorem]{Definition}

\theoremstyle{remark}
\newtheorem{remark}[theorem]{Remark}

\numberwithin{equation}{section}

\begin{document}

\maketitle

\begin{abstract}Given a  generalized $\overline M^{n+1}=I\times_{\phi}F^n$
Robertson-Walker spacetime  we will classify strongly stable
spacelike hypersurfaces  with constant mean curvature whose warping
function verifies a certain convexity condition. More precisely, we
will show that given $x:M^n\rightarrow\overline M^{n+1}$ a closed
spacelike hypersurfaces of $\overline M^{n+1}$ with constant mean
curvature $H$ and the warping function $\phi$ satisfying
$\phi''\geq\max\{H\phi',0\}$, then $M^{n}$ is either minimal or a
spacelike slice $M_{t_0}=\{t_0\}\times F$, for some $t_0\in I$.

\end{abstract}

\section{Introduction}

Spacelike hypersurfaces with constant mean curvature in Lorentz
manifolds have been object of great interest in recent years,
both from physical and mathematical points of view. In~\cite{ABC:03},
the authors studied the uniqueness of spacelike hypersurfaces  with CMC in
generalized Robertson-Walker (GRW) spacetimes, namely, Lorentz warped
products with 1-dimensional negative  definite base and Riemannian fiber.
They proved that in a GRW spacetime obeying the timelike convergence condition
(i.e, the Ricci curvature is non-negative on timelike directions), every compact
spacelike hypersurface with CMC must be umbilical. Recently, Al\'{i}as and  Montiel obtained,
in~\cite{AM:01}, a more general condition on the warping function $f$ that is
sufficient in order to guarantee uniqueness. More precisely, they proved the following

\begin{theorem} Let $f:I\rightarrow\mathbb R$ be a positive smooth function defined on
an open interval, such that $ff^{''}-(f^{'})^{2}\leq 0$, that is,
such that $-\log f$ is convex. Then, the only compact spacelike
hypersurfaces immersed into a generalized Robertson-Walker spacetime
$I\times_fF^{n}$ and having constant mean curvature are the slices
$\{t\}\times F$, for a {\rm(}necessarily compact{\rm)} Riemannian
manifold $F$.
\end{theorem}

Stability questions concerning CMC, compact hypersurfaces in
Riemannian space forms began with Barbosa and do Carmo
in~\cite{BdC:84}, and Barbosa, Do Carmo and Eschenburg
in~\cite{BdCE:88}. In the former paper, they introduced the notion
of stability and proved that spheres are the only stable critical
points for the area functional, for volume-preserving variations. In
the setting of spacelike hypersurfaces in Lorentz manifolds, Barbosa
and Oliker proved in~\cite{Barbosa:93} that CMC spacelike
hypersurfaces are critical points of volume-preserving variations.
Moreover, by computing the second variation formula they showed that
CMC embedded spheres in the de Sitter space $S_1^{n+1}$ maximize the
area functional for such variations. In this paper, we give a
characterization of {\em strongly stable}, CMC spacelike
hypersurfaces in GRW spacetimes, the essential tool for the proof
being a formula for the Laplacian of a new support function. More
precisely, it is our purpose to show the following

\begin{theorem}
Let $\overline M^{n+1}=I\times_{\phi}F^n$ be a generalized
Robertson-Walker spacetime, and $x:M^n\rightarrow\overline M^{n+1}$
be a closed spacelike hypersurface of $\overline M^{n+1}$, having
constant mean curvature $H$. If the warping function $\phi$
satisfies $\phi''\geq\max\{H\phi',0\}$ and $M^n$ is strongly
stable, then $M^{n}$ is either minimal or a spacelike slice
$M_{t_0}=\{t_0\}\times F$, for some $t_0\in I$.
\end{theorem}

\section{Stable spacelike hypersurfaces}

In what follows, $\overline M^{n+1}$ denotes an orientable, time-oriented Lorentz manifold
with Lorentz metric $\overline g=\langle\,\,,\,\,\rangle$ and semi-Riemannian connection
$\overline\nabla$. If $x:M^n\rightarrow\overline M^{n+1}$ is a spacelike hypersurface of
$\overline M^{n+1}$, then $M^n$ is automatically orientable (\cite{O'Neill:83}, p. 189),
and one can choose a globally defined unit normal vector field $N$ on $M^n$ having the
same time-orientation of $V$, that is, such that
$$\langle V,N\rangle<0$$
on $M$. One says that such an $N$ {\em points to the future}.

A {\em variation} of $x$ is a smooth map
$$X:M^n\times(-\epsilon,\epsilon)\rightarrow\overline M^{n+1}$$
satisfying the following conditions:
\begin{enumerate}
\item[(1)] For $t\in(-\epsilon,\epsilon)$, the map $X_t:M^n\rightarrow\overline M^{n+1}$
given by $X_t(p)=X(t,p)$ is a spacelike immersion such that $X_0=x$.
\item[(2)] $X_t\big|_{\partial M}=x\big|_{\partial M}$, for all $t\in(-\epsilon,\epsilon)$.
\end{enumerate}

The {\em variational field} associated to the variation $X$ is the vector field
$\frac{\partial X}{\partial t}$. Letting
$f=-\langle\frac{\partial X}{\partial t},N\rangle$, we get
$$\frac{\partial X}{\partial t}\Big|_M=fN+\left(\frac{\partial X}{\partial t}\right)^T,$$
where $T$ stands for tangential components. The {\em balance of volume} of the variation $X$
is the function $\mathcal V:(-\epsilon,\epsilon)\rightarrow\mathbb R$ given by
$$\mathcal V(t)=\int_{M\times[0,t]}X^*(d\overline M),$$
where $d\overline M$ denotes the volume element of $\overline M$.

The {\em area functional} $\mathcal A:(-\epsilon,\epsilon)\rightarrow\mathbb R$ associated
to the variation $X$ is given by
$$\mathcal A(t)=\int_MdM_t,$$
where $dM_t$ denotes the volume element of the metric induced in $M$ by $X_t$. Note that
$dM_0=dM$ and $\mathcal A(0)=\mathcal A$, the volume of $M$.
The following lemma is classical:

\begin{lemma}\label{lemma:first variation}
Let $\overline M^{n+1}$ be a time-oriented Lorentz manifold and
$x:M^n\rightarrow\overline M^{n+1}$ a spacelike closed hypersurface having mean curvature $H$.
If $X:M^n\times(-\epsilon,\epsilon)\rightarrow\overline M^{n+1}$ is a variation of $x$, then
$$\frac{d\mathcal V}{dt}\Big|_{t=0}=\int_MfdM,\ \   \ \ \frac{d\mathcal A}{dt}\Big|_{t=0}=\int_MnHfdM.$$
\end{lemma}

Set $H_0=\frac{1}{\mathcal A}\int_MdM$ and
$\mathcal J:(-\epsilon,\epsilon)\rightarrow\mathbb R$ given by
$$\mathcal J(t)=\mathcal A(t)-nH_0\mathcal V(t).$$
$\mathcal J$ is called the {\em Jacobi functional} associated to the variation,
and it is a well known result~\cite{BdCE:88} that $x$ has constant mean curvature $H_0$
if and only if $\mathcal J'(0)=0$ for all variations $X$ of $x$.

We wish to study here immersions $x:M^n\rightarrow\overline M^{n+1}$ that
maximize $\mathcal J$ for all variations $X$. Since $x$ must be a critical point of
$\mathcal J$, it thus follows from the above discussion that $x$ must have constant mean
curvature. Therefore, in order to examine whether or not some critical immersion $x$ is
actually a maximum for $\mathcal J$, one certainly needs to study the second variation
$\mathcal J''(0)$. We start with the following

\begin{proposition} Let $x:M^n\rightarrow\overline M^{n+1}$ be a closed spacelike
hypersurface of the time-oriented Lorentz manifold $\overline M^{n+1}$, and
$X:M^n\times(-\epsilon,\epsilon)\rightarrow\overline M^{n+1}$ be a variation of $x$. Then,
\begin{equation}\label{eq:fundamental relation}
n\frac{\partial H}{\partial t}=\Delta f-\left\{\overline{Ric}(N,N)+|A|^2\right\}f-n\langle\left(\frac{\partial X}{\partial t}\right)^T,\nabla H\rangle.
\end{equation}
\end{proposition}

Although the above proposition is known to be true, we believe there is a lack,
in the literature, of a clear proof of it in this degree of generality, so we present
a simple proof here.

\begin{proof} Let $p\in M$ and $\{e_k\}$ be a moving frame on a neighborhood $U\subset M$ of
$p$, geodesic at $p$ and diagonalizing $A$ at $p$, with $Ae_k=\lambda_ke_k$ for
$1\leq k\leq n$. Extend $N$ and the $e_k's$ to a neighborhood of $p$ in $\overline M$, so
that $\langle N,e_k\rangle=0$ and $(\overline\nabla_Ne_k)(p)=0$. Then
\begin{eqnarray*}
n\frac{\partial H}{\partial t}&=&-{\rm tr}\left(\frac{\partial A}{\partial t}\right)=-\sum_k\langle\frac{\partial A}{\partial t}e_k,e_k\rangle=-\sum_k\langle\left(\overline\nabla_{\frac{\partial X}{\partial t}}A\right)e_k,e_k\rangle\\
&=&-\sum_k\left\{\langle\overline\nabla_{\frac{\partial X}{\partial t}}Ae_k,e_k\rangle-\langle A\overline\nabla_{\frac{\partial X}{\partial t}}e_k,e_k\rangle\right\}\\
&=&\sum_k\langle\overline\nabla_{\frac{\partial X}{\partial t}}\overline\nabla_{e_k}N,e_k\rangle+\sum_k\langle A\overline\nabla_{e_k}\frac{\partial X}{\partial t},e_k\rangle,
\end{eqnarray*}
where in the last equality we used the fact that $[\frac{\partial X}{\partial t},e_k]=0$.
Letting
$$I=\sum_k\langle\overline\nabla_{\frac{\partial X}{\partial t}}\overline\nabla_{e_k}N,e_k\rangle\ \ \text{and}\ \ II=\sum_k\langle A\overline\nabla_{e_k}\frac{\partial X}{\partial t},e_k\rangle,$$
we have
\begin{eqnarray*}
I&=&\sum_k\left\{\langle\overline\nabla_{\frac{\partial X}{\partial t}}\overline\nabla_{e_k}N-\overline\nabla_{e_k}\overline\nabla_{\frac{\partial X}{\partial t}}N+\overline\nabla_{[e_k,\frac{\partial X}{\partial t}]}N,e_k\rangle+\overline\nabla_{e_k}\overline\nabla_{\frac{\partial X}{\partial t}}N,e_k\rangle\right\}\\
&=&\sum_k\left\{\langle\overline R\left(e_k,\frac{\partial X}{\partial t}\right)N,e_k\rangle+\langle\overline\nabla_{e_k}\overline\nabla_{\frac{\partial X}{\partial t}}N,e_k\rangle\right\}\\
&=&-\overline{Ric}\left(\frac{\partial X}{\partial t},N\right)+\sum_k\langle\overline\nabla_{e_k}\overline\nabla_{\frac{\partial X}{\partial t}}N,e_k\rangle.
\end{eqnarray*}

Since the frame $\{e_k\}$ is geodesic at $p$, it follows that
$$\langle\overline\nabla_{\frac{\partial X}{\partial t}}N,\overline\nabla_{e_k}e_k\rangle=\langle\overline\nabla_{\frac{\partial X}{\partial t}}N,N\rangle\langle\overline\nabla_{e_k}e_k,N\rangle=0$$
at $p$, and hence
\begin{eqnarray*}
\langle\overline\nabla_{e_k}\overline\nabla_{\frac{\partial X}{\partial t}}N,e_k\rangle&=&e_k\langle\overline\nabla_{\frac{\partial X}{\partial t}}N,e_k\rangle=-e_k\langle N,\overline\nabla_{\frac{\partial X}{\partial t}}e_k\rangle=-e_k\langle N,\overline\nabla_{e_k}\frac{\partial X}{\partial t}\rangle\\
&=&-e_ke_k\langle N,\frac{\partial X}{\partial t}\rangle+e_k\langle\overline\nabla_{e_k}N,\frac{\partial X}{\partial t}\rangle\\
&=&e_ke_k(f)+e_k\langle\overline\nabla_{e_k}N,\left(\frac{\partial X}{\partial t}\right)^T\rangle\\
&=&e_ke_k(f)+\langle\overline\nabla_{e_k}\overline\nabla_{e_k}N,\left(\frac{\partial X}{\partial t}\right)^T\rangle-\langle Ae_k,\overline\nabla_{e_k}\left(\frac{\partial X}{\partial t}\right)^T\rangle.
\end{eqnarray*}

For $II$, we have
\begin{eqnarray*}
II&=&\sum_k\langle Ae_k,\overline\nabla_{e_k}\frac{\partial X}{\partial t}\rangle=\sum_k\langle Ae_k,\overline\nabla_{e_k}(fN+\left(\frac{\partial X}{\partial t}\right)^T)\rangle\\
&=&\sum_k\langle Ae_k,f\overline\nabla_{e_k}N\rangle+\sum_k\langle Ae_k,\overline\nabla_{e_k}\left(\frac{\partial X}{\partial t}\right)^T\rangle\\
&=&-f|A|^2+\sum_k\langle Ae_k,\overline\nabla_{e_k}\left(\frac{\partial X}{\partial t}\right)^T\rangle\\
\end{eqnarray*}

Therefore,
\begin{equation}\label{eq:aux_I}
n\frac{\partial H}{\partial t}=-\overline{Ric}\left(\frac{\partial X}{\partial t},N\right)+\Delta f-f|A|^2+\sum_k\langle\overline\nabla_{e_k}\overline\nabla_{e_k}N,\left(\frac{\partial X}{\partial t}\right)^T\rangle.
\end{equation}

Now, letting
$$\frac{\partial X}{\partial t}=\sum_l^n\alpha_le_l+fN$$
and $Ae_k=\sum_jh_{jk}e_j$, one successively gets
\begin{eqnarray*}
\overline{Ric}\left(\frac{\partial X}{\partial t},N\right)&=&\sum_l\alpha_l\overline{Ric}(N,e_l)+f\overline{Ric}(N,N)\\
&=&\sum_{k,l}\alpha_l\langle\overline R(e_k,e_l)e_k,N\rangle+f\overline{Ric}(N,N)\\
\end{eqnarray*}
and, since $(\overline\nabla_Ne_k)(p)=0$,
\begin{eqnarray*}
\langle\overline R(e_k,e_l)e_k,N\rangle_p&=&\langle\overline\nabla_{e_l}\overline\nabla_{e_k}e_k-\overline\nabla_{e_k}\overline\nabla_{e_l}e_k,N\rangle_p\\
&=&e_l\langle\overline\nabla_{e_k}e_k,N\rangle_p-\langle\overline\nabla_{e_k}e_k,\overline\nabla_{e_l}N\rangle_p-e_k\langle\overline\nabla_{e_l}e_k,N\rangle_p\\
&=&-e_l\langle e_k,\overline\nabla_{e_k}N\rangle_p+e_k\langle e_k,\overline\nabla_{e_l}N\rangle_p\\
&=&e_l(h_{kk})-e_k(h_{kl}),
\end{eqnarray*}
so that
\begin{equation}\label{eq:aux_II}
\overline{Ric}\left(\frac{\partial X}{\partial t},N\right)_p=\sum_{k,l}\alpha_le_l(h_{kk})-\sum_{k,l}\alpha_le_k(h_{kl})+f\overline{Ric}(N,N)_p.
\end{equation}

Also,
\begin{eqnarray*}
\alpha_l\langle\overline\nabla_{e_k}\overline\nabla_{e_k}N,e_l\rangle&=&\alpha_l\langle\nabla_{e_k}\overline\nabla_{e_k}N,e_l\rangle=-\alpha_l\sum_j\langle\nabla_{e_k}h_{kj}e_j,e_l\rangle\\
&=&-\alpha_l\sum_j\left\{e_k(h_{kj})\delta_{lj}+h_{kj}\langle\nabla_{e_k}e_j,e_l\rangle\right\}\\
&=&-\alpha_le_k(h_{kl}),
\end{eqnarray*}
and hence
\begin{equation}\label{eq:aux_III}
\sum_k\langle\overline\nabla_{e_k}\overline\nabla_{e_k}N,\left(\frac{\partial X}{\partial t}\right)^T\rangle=-\sum_{k,l}\alpha_le_k(h_{kl}).
\end{equation}

Substituting (\ref{eq:aux_II}) and (\ref{eq:aux_III}) into (\ref{eq:aux_I}), we finally arrive
at
\begin{eqnarray*}
n\frac{\partial H}{\partial t}&=&-\sum_{k,l}\alpha_le_l(h_{kk})-f\overline{Ric}(N,N)_p+\Delta f-f|A|^2\\
&=&-\left(\frac{\partial X}{\partial t}\right)^T(nH)-f\overline{Ric}(N,N)_p+\Delta f-f|A|^2.
\end{eqnarray*}
\end{proof}

\begin{proposition} Let $\overline M^{n+1}$ be a Lorentz manifold and
$x:M^n\rightarrow\overline M^{n+1}$ be a closed spacelike hypersurface having constant
mean curvature $H$. If $X:M^n\times(-\epsilon,\epsilon)\rightarrow\overline M^{n+1}$
is a variation of $x$, then
\begin{equation}\label{eq:second variation of J}
\mathcal J''(0)(f)=\int_Mf\left\{\Delta f-\left(\overline{Ric}(N,N)+|A|^2\right)f\right\}dM.
\end{equation}
\end{proposition}

\begin{proof} In the notations of the above discussion, set $f=f(0)$ and note that $H(0)=H$.
It follows from lemma~\ref{lemma:first variation} that
$$\mathcal J'(t)=\int_Mn\left\{H(t)-H\right\}f(t)dM_t.$$
Therefore, differentiating with respect to $t$ once more
\begin{eqnarray*}
\mathcal J''(0)&=&\int_MnH'(0)f(0)dM_0+\int_Mn\left\{H(0)-H\right\}\frac{d}{dt}f(t)dM_t\Big|_{t=0}\\
&=&\int_MnH'(0)fdM.
\end{eqnarray*}

Taking into account that $H$ is constant, relation (\ref{eq:fundamental relation}) finally
gives formula~\ref{eq:second variation of J}
\end{proof}

It follows from the previous result that $\mathcal J''(0)=\mathcal J''(0)(f)$ depends only
on $f\in C^{\infty}(M)$, for which there exists a variation $X$ of $M^n$ such that
$\left(\frac{\partial X}{\partial t}\right)^{\bot}=fN$. Therefore, the following definition
makes sense:

\begin{definition} Let $\overline M^{n+1}$ be a Lorentz manifold and
$x:M^n\rightarrow\overline M^{n+1}$ be a closed spacelike hypersurface having constant
mean curvature $H$. We say that $x$ is strongly stable if, for every function
$f\in C^{\infty}(M)$ for which there exists a variation $X$ of $M^n$ such that
$\left(\frac{\partial X}{\partial t}\right)^{\bot}=fN$, one has $\mathcal J''(0)(f)\leq 0$.
\end{definition}

\section{Conformal vector fields}

As in the previous section, let $\overline M^{n+1}$ be a Lorentz manifold.
A vector field $V$ on $\overline M^{n+1}$ is said to be {\em conformal} if
\begin{equation}
\mathcal L_V\langle\,\,,\,\,\rangle=2\psi\langle\,\,,\,\,\rangle
\end{equation}
for some function $\psi\in C^{\infty}(\overline M)$, where $\mathcal L$ stands for the Lie
derivative of the Lorentz metric of $\overline M$. The function $\psi$ is called the
{\em conformal factor} of $V$.

Since $\mathcal L_V(X)=[V,X]$ for all $X\in\mathcal X(\overline M)$, it follows from the
tensorial character of $\mathcal L_V$ that $V\in\mathcal X(\overline M)$ is conformal if
and only if
\begin{equation}\label{eq:1.1}
\langle\overline\nabla_XV,Y \rangle+\langle X,\overline\nabla_YV\rangle=2\psi\langle X,Y\rangle,
\end{equation}
for all $X,Y\in\mathcal X(\overline M)$. In particular, $V$ is a Killing vector field
relatively to $\overline g$ if and only if $\psi\equiv 0$.

Any Lorentz manifold $\overline M^{n+1}$, possessing a globally defined, timelike conformal
vector field is said to be a {\em conformally stationary spacetime}.

\begin{proposition}\label{prop:Laplacian of conformal vector field}
Let $\overline M^{n+1}$ be a conformally stationary Lorentz manifold, with conformal vector
field $V$ having conformal factor $\psi:\overline M^{n+1}\rightarrow\mathbb R$. Let also
$x:M^n\rightarrow\overline M^{n+1}$ be a spacelike hypersurface of $\overline M^{n+1}$,
and $N$ a future-pointing, unit normal vector field globally defined on $M^n$.
If $f=\langle V,N\rangle$, then
\begin{equation}\label{eq:Laplacian formula_I}
\Delta f=n\langle V,\nabla H\rangle+f\left\{\overline{Ric}(N,N)+|A|^2\right\}+n\left\{H\psi-N(\psi)\right\},
\end{equation}
where $\overline{Ric}$ denotes the Ricci tensor of $\overline M$, $A$ is the second fundamental
form of $x$ with respect to $N$, $H=-\frac{1}{n}{\rm tr}(A)$ is the mean curvature of $x$ and
$\nabla H$ denotes the gradient of $H$ in the metric of $M$.
\end{proposition}

\begin{proof} Fix $p\in M$ and let $\{e_k\}$ be an orthonormal moving frame on $M$, geodesic
at $p$. Extend the $e_k$ to a neighborhood of $p$ in $\overline M$, so that
$(\overline\nabla_Ne_k)(p)=0$, and let
$$V=\sum_l^n\alpha_le_l-fN.$$
Then
\begin{eqnarray*}
f=\langle N,V\rangle\Rightarrow e_k(f)&=&\langle\overline\nabla_{e_k}N,V\rangle+\langle N,\overline\nabla_{e_k}V\rangle\\
&=&-\langle Ae_k,V\rangle+\langle N,\overline\nabla_{e_k}V\rangle,
\end{eqnarray*}
so that
\begin{eqnarray}\label{eq:I}
\Delta f&=&\sum_ke_k(e_k(f))=-\sum_ke_k\langle Ae_k,V\rangle+\sum_ke_k\langle N,\overline\nabla_{e_k}V\rangle\nonumber\\
&=&-\sum_k\langle\overline\nabla_{e_k}Ae_k,V\rangle-2\sum_k\langle Ae_k,\overline\nabla_{e_k}V\rangle+\sum_k\langle N,\overline\nabla_{e_k}\overline\nabla_{e_k}V\rangle.
\end{eqnarray}

Now, differentiating $Ae_k=\sum_lh_{kl}e_l$ with respect to $e_k$, one gets at $p$
\begin{eqnarray}\label{eq:II}
\sum_k\langle\overline\nabla_{e_k}Ae_k,V\rangle&=&\sum_{k,l}e_k(h_{kl})\langle e_l,V\rangle+\sum_{k,l}h_{kl}\langle\overline\nabla_{e_k}e_l,V\rangle\nonumber\\
&=&\sum_{k,l}\alpha_le_k(h_{kl})-\sum_{k,l}h_{kl}\langle\overline\nabla_{e_k}e_l,N\rangle\langle V,N\rangle\nonumber\\
&=&\sum_{k,l}\alpha_le_k(h_{kl})-\sum_{k,l}h_{kl}^2f\nonumber\\
&=&\sum_{k,l}\alpha_le_k(h_{kl})-f|A|^2.
\end{eqnarray}

Asking further that $Ae_k=\lambda_ke_k$ at $p$ (which is always possible), we have at $p$
\begin{equation}\label{eq:III}
\sum_k\langle Ae_k,\overline\nabla_{e_k}V\rangle=\sum_k\lambda_k\langle e_k,\overline\nabla_{e_k}V\rangle=\sum_k\lambda_k\psi=-nH\psi.
\end{equation}

In order to compute the last summand of (\ref{eq:I}), note that the conformality of $V$ gives
$$\langle\overline\nabla_NV,e_k\rangle+\langle N,\overline\nabla_{e_k}V\rangle=0$$
for all $k$. Hence, differentiating the above relation in the direction of $e_k$, we get
$$\langle\overline\nabla_{e_k}\overline\nabla_NV,e_k\rangle+\langle\overline\nabla_NV,\overline\nabla_{e_k}e_k\rangle+\langle\overline\nabla_{e_k}N,\overline\nabla_{e_k}V\rangle+\langle N,\overline\nabla_{e_k}\overline\nabla_{e_k}V\rangle=0.$$
However, at $p$ one has
\begin{eqnarray*}
\langle\overline\nabla_NV,\overline\nabla_{e_k}e_k\rangle&=&-\langle\overline\nabla_NV,\langle\overline\nabla_{e_k}e_k,N\rangle N\rangle=-\langle\overline\nabla_NV,\lambda_kN\rangle\\
&=&-\lambda_k\psi\langle N,N\rangle=\lambda_k\psi
\end{eqnarray*}
and
$$\langle\overline\nabla_{e_k}N,\overline\nabla_{e_k}V\rangle=-\lambda_k\langle e_k,\overline\nabla_{e_k}V\rangle=-\lambda_k\psi,$$
so that
\begin{equation}\label{eq:IV}
\langle\overline\nabla_{e_k}\overline\nabla_NV,e_k\rangle+\langle N,\overline\nabla_{e_k}\overline\nabla_{e_k}V\rangle=0
\end{equation}
at $p$. On the other hand, since
$$[N,e_k](p)=(\overline\nabla_Ne_k)(p)-(\overline\nabla_{e_k}N)(p)=\lambda_ke_k(p),$$
it follows from (\ref{eq:IV}) that
\begin{eqnarray*}
\langle\overline R(N,e_k)V,e_k\rangle_p&=&\langle\overline\nabla_{e_k}\overline\nabla_NV-\overline\nabla_N\overline\nabla_{e_k}V+\overline\nabla_{[N,e_k]}V,e_k\rangle_p\\
&=&-\langle N,\overline\nabla_{e_k}\overline\nabla_{e_k}V\rangle_p-N\langle\overline\nabla_{e_k}V,e_k\rangle_p+\langle\overline\nabla_{\lambda_ke_k}V,e_k\rangle_p\\
&=&-\langle N,\overline\nabla_{e_k}\overline\nabla_{e_k}V\rangle_p-N(\psi)+\lambda_k\psi,
\end{eqnarray*}
and hence
\begin{equation}\label{eq:V}
\sum_k\langle N,\overline\nabla_{e_k}\overline\nabla_{e_k}V\rangle_p=-nN(\psi)-nH\psi-\overline{Ric}(N,V)_p
\end{equation}
Finally,
\begin{eqnarray*}
\overline{Ric}(N,V)&=&\sum_l\alpha_l\overline{Ric}(N,e_l)-f\overline{Ric}(N,N)\\
&=&\sum_{k,l}\alpha_l\langle\overline R(e_k,e_l)e_k,N\rangle-f\overline{Ric}(N,N),
\end{eqnarray*}
and
\begin{eqnarray*}
\langle\overline R(e_k,e_l)e_k,N\rangle_p&=&\langle\overline\nabla_{e_l}\overline\nabla_{e_k}e_k-\overline\nabla_{e_k}\overline\nabla_{e_l}e_k,N\rangle_p\\
&=&e_l\langle\overline\nabla_{e_k}e_k,N\rangle_p-\langle\overline\nabla_{e_k}e_k,\overline\nabla_{e_l}N\rangle_p-e_k\langle\overline\nabla_{e_l}e_k,N\rangle_p\\
&&+\langle\overline\nabla_{e_l}e_k,\overline\nabla_{e_k}N\rangle_p\\
&=&-e_l\langle e_k,\overline\nabla_{e_k}N\rangle_p+e_k\langle e_k,\overline\nabla_{e_l}N\rangle_p\\
&=&e_l(h_{kk})-e_k(h_{kl}),
\end{eqnarray*}
so that
$$\overline{Ric}(N,V)_p=\sum_{k,l}\alpha_le_l(h_{kk})-\sum_{k,l}\alpha_le_k(h_{kl})-f\overline{Ric}(N,N)_p,$$
and it follows from (\ref{eq:V}) that
\begin{eqnarray}\label{eq:VI}
\sum_k\langle N,\overline\nabla_{e_k}\overline\nabla_{e_k}V\rangle_p&=&-nN(\psi)-nH\psi+V^T(nH)\nonumber\\
&&+\sum_{k,l}\alpha_le_k(h_{kl})+f\overline{Ric}(N,N).
\end{eqnarray}

Substituting (\ref{eq:II}), (\ref{eq:III}) and (\ref{eq:VI}) into (\ref{eq:I}), one gets
the desired formula (\ref{eq:Laplacian formula_I}).
\end{proof}

\section{Applications}

A particular class of conformally stationary spacetimes is that of
{\em generalized Robertson-Walker} spacetimes~\cite{ABC:03}, namely, warped products
$\overline M^{n+1}=I\times_{\phi}F^n$, where $I\subseteq\mathbb R$ is an interval with the
metric $-dt^2$, $F^n$ is an $n$-dimensional Riemannian manifold and
$\phi:I\rightarrow\mathbb R$ is positive and smooth. For such a space, let
$\pi_I:\overline M^{n+1}\rightarrow I$ denote the canonical projection onto the $I-$factor.
Then the vector field
$$V=(\phi\circ\pi_I)\frac{\partial}{\partial t}$$
is conformal, timelike and closed (in the sense that its dual
$1-$form is closed), with conformal factor $\psi=\phi'$, where the
prime denotes differentiation with respect to $t$.
Moreover, according to~\cite{Montiel:99}, for $t_0\in I$, orienting the
(spacelike) leaf $M_{t_0}^n=\{t_0\}\times F^n$ by using the
future-pointing unit normal vector field $N$, it follows that
$M_{t_0}$ has constant mean curvature
$$H=\frac{\phi'(t_0)}{\phi(t_0)}.$$

If $\overline M^{n+1}=I\times_{\phi}F^n$ is a generalized Robertson-Walker spacetime and
$x:M^n\rightarrow\overline M^{n+1}$ is a complete spacelike hypersurface of
$\overline M^{n+1}$, such that $\phi\circ\pi_I$ is limited on $M$, then
$\pi_F\big|_M:M^n\rightarrow F^n$ is necessarily a covering map (\cite{ABC:03}).
In particular, if $M^n$ is closed, then $F^n$ is automatically closed.

One has the following corollary of proposition~\ref{prop:Laplacian of conformal vector field}:

\begin{corollary}
Let $\overline M^{n+1}=I\times_{\phi}F^n$ be a generalized Robertson-Walker spacetime,
and $x:M^n\rightarrow\overline M^{n+1}$ a spacelike hypersurface of $\overline M^{n+1}$,
having constant mean curvature $H$. Let also $N$ be a future-pointing unit normal vector
field globally defined on $M^n$. If $V=(\phi\circ\pi_I)\frac{\partial}{\partial t}$ and
$f=\langle V,N\rangle$, then
\begin{equation}\label{eq:Laplacian formula_II}
\Delta f=\left\{\overline{Ric}(N,N)+|A|^2\right\}f+n\left\{H\phi'+\phi''\langle N,\frac{\partial}{\partial t}\rangle\right\}.
\end{equation}
where $\overline{Ric}$ denotes the Ricci tensor of $\overline M$, $A$ is the second fundamental
form of $x$ with respect to $N$, and $H=-\frac{1}{n}{\rm tr}(A)$ is the mean curvature of $x$.
\end{corollary}

\begin{proof} First of all,
$f=\langle V,N\rangle=\phi\langle N,\frac{\partial}{\partial t}\rangle$, and
it thus follows from (\ref{eq:Laplacian formula_I}) that
$$\Delta f=\left\{\overline{Ric}(N,N)+|A|^2\right\}f+n\left\{H\phi'-N(\phi')\right\}.$$
However,
$$\overline\nabla\phi'=-\langle\overline\nabla\phi',\frac{\partial}{\partial t}\rangle\frac{\partial}{\partial t}=-\phi''\frac{\partial}{\partial t},$$
so that
$$N(\phi')=\langle N,\overline\nabla\phi'\rangle=-\phi''\langle N,\frac{\partial}{\partial t}\rangle$$
\end{proof}

We can now state and prove our main result:

\begin{theorem}
Let $\overline M^{n+1}=I\times_{\phi}F^n$ be a generalized
Robertson-Walker spacetime, and $x:M^n\rightarrow\overline M^{n+1}$
be a closed spacelike hypersurface of $\overline M^{n+1}$, having
constant mean curvature $H$. If the warping function $\phi$
satisfies $\phi''\geq\max\{H\phi',0\}$ and $M^n$ is strongly
stable, then $M^{n}$ is either minimal or a spacelike slice
$M_{t_0}=\{t_0\}\times F$, for some $t_0\in I$.
\end{theorem}

\begin{proof} Since $M^n$ is strongly stable, we have
$$0\geq\mathcal J''(0)(g)=\int_Mg\left\{\Delta g-\left(\overline{Ric}(N,N)+|A|^2\right)g\right\}dM$$
for all $g\in C^{\infty}(M)$ for which $gN$ is the normal component of the variational field of
some variation of $M^n$. In particular, if
$f=\langle V,N\rangle=\phi\langle N,\frac{\partial}{\partial t}\rangle$, where
$V=(\phi\circ\pi_I)\frac{\partial}{\partial t}$, and $g=-f=-\langle V,N\rangle$, then
$$\Delta g=\left\{\overline{Ric}(N,N)+|A|^2\right\}g-n\left\{H\phi'+\phi''\langle N,\frac{\partial}{\partial t}\rangle\right\}.$$
Therefore, $M^n$ stable implies
$$0\geq\int_M\phi\langle N,\frac{\partial}{\partial t}\rangle\left\{H\phi'+\phi''\langle N,\frac{\partial}{\partial t}\rangle\right\}dM$$
Letting $\theta$ be the hyperbolic angle between $N$ and $\frac{\partial}{\partial t}$, it
follows from the reversed Cauchy-Schwarz inequality that
$\cosh\theta=-\langle N,\frac{\partial}{\partial t}\rangle$, with $\cosh\theta\equiv 1$ if and only
if $N$ and $\frac{\partial}{\partial t}$ are collinear at every point, that is, if and only if
$M^n$ is a spacelike leaf $M_{t_0}$ for some $t_0\in I$. Hence,
$$0\geq\int_M\phi\cosh\theta\left\{-H\phi'+\phi''\cosh\theta\right\}dM.$$

Now, notice that $-H\phi'+\phi''\cosh\theta\geq-\phi''+\phi''\cosh\theta$, which gives
$$\phi\cosh\theta(-H\phi'+\phi''\cosh\theta)\geq\phi\phi''\cosh\theta(\cosh\theta-1).$$
Therefore,
$$0\geq\int_M\phi\cosh\theta(-H\phi'+\phi''\cosh\theta)dM\geq\int_M\phi\phi''\cosh\theta(\cosh\theta-1)\geq 0,$$
and hence
$$\phi''(\cosh\theta-1)=0\ \ \text{and}\ \ \phi''=H\phi'$$
on $M$. If, for some $p\in M$, one has $\phi''(p)=0$, then $\phi'H=0$ at $p$. If $H\neq 0$,
then $\phi'(p)=0$. But if this is the case, then
proposition 7.35 of~\cite{O'Neill:83} gives that
$$\overline\nabla_V\frac{\partial}{\partial t}=\frac{\phi'}{\phi}V=0$$
at $p$ for any $V$, and $M$ is totally geodesic at $p$. In particular, $H=0$, a
contradiction. Therefore, either $\phi''(p)=0$ for some $p\in M$, and $M$ is minimal,
or $\phi''\neq 0$ on all of $M$, whence $\cosh\theta=1$ always, and $M$ is an umbilical
leaf such that $\phi''=H\phi'$.
\end{proof}

\remark Note that $\frac{\phi''}{\phi'}=H=\frac{\phi'}{\phi}$, i.e.,
$\phi''\phi-(\phi')^2=0$, which is a limit case of Al\'{i}as and Montiel's timelike convergent condition.


\begin{thebibliography}{10}

\bibitem{ABC:03} L. J. Al\'{\i}as, A. Brasil Jr. and A. G. Colares,
{\em Integral Formulae for Spacelike Hypersurfaces in Conformally Stationary Spacetimes and Applications},
Proc. Edinburgh Math. Soc. 46, (2003) 465-488 .

\bibitem{AM:01} L. J. Al\'{i}as \and S. Montiel,
{\em Uniqueness of Spacelike Hypersurfaces with Constant Mean Curvature in Generalized Robertson-Walker Spacetimes},
Proceedings of the International Conference held to honour the 60th birthday of A.M.Naveira, World Scientific, (2001) 59-69.

\bibitem{Barbosa:97} J. L. M. Barbosa \and A. G. Colares,
{\em Stability of Hypersurfaces with Constant $r-$Mean Curvature},
Ann. Global Anal. Geom. 15, (1997) 277-297.

\bibitem{BdC:84} J. L. M. Barbosa \and M. do Carmo,
{\em Stability of Hypersurfaces with Constant Mean Curvature},
Math. Z. 185, (1984) 339-353.

\bibitem{BdCE:88} J. L. M. Barbosa, M. do Carmo \and J. Eschenburg,
{\em Stability of Hypersurfaces with Constant Mean Curvature},
Math. Z. 197, (1988) 123-138.

\bibitem{Barbosa:93} J. L. M. Barbosa \and V. Oliker,
{\em Spacelike Hypersurfaces with Constant Mean Curvature in Lorentz Spaces},
Matem. Contempor\^anea 4, (1993) 27-44.

\bibitem{Montiel:99} S. Montiel, {\em Uniqueness of Spacelike Hypersurfaces of Constant Mean 
Curvature in foliated Spacetimes}, Math. Ann. 314, (1999) 529-553. 


\bibitem{O'Neill:83} B. O'Neill,
{\em Semi-Riemannian Geometry with Applications to Relativity},
London, Academic Press (1983).
\end{thebibliography}
\end{document}